\documentclass{amsart}
\usepackage{amsmath}
\usepackage{latexsym,ulem}
\usepackage{amssymb}
\usepackage{graphics}

 \input epsf%
\newtheorem{theorem}{Theorem}[section]

\newtheorem{lemma}{Lemma}[section]
\newtheorem{remark}{Remark}[section]

\newtheorem{Definition}{Definition}[section]
\def\qed{{\hfill{\vrule height4pt width3pt depth2pt}}}

   \long\def\comment#1{}

\def\ad#1{\begin{aligned}#1\end{aligned}}  \def\b#1{\mathbf{#1}}
\def\a#1{\begin{align*}#1\end{align*}} \def\an#1{\begin{align}#1\end{align}}
\def\e#1{\begin{equation}#1\end{equation}} \def\d{\operatorname{div}}
\def\p#1{\begin{pmatrix}#1\end{pmatrix}} 
  \numberwithin{equation}{section}
\numberwithin{table}{section}
\numberwithin{figure}{section}

\def\boxit#1{\vbox{\hrule height1pt \hbox{\vrule width1pt\kern1pt
     #1\kern1pt\vrule width1pt}\hrule height1pt }}
 \def\lab#1{\boxit{\small #1}\label{#1}}
  \def\mref#1{\boxit{\small #1}\ref{#1}}
 \def\meqref#1{\boxit{\small #1}\eqref{#1}}

  \def\lab#1{\label{#1}} \def\mref#1{\ref{#1}} \def\meqref#1{\eqref{#1}}

\begin{document}

\newcommand{\disp}{\displaystyle}
\newcommand{\eps}{\varepsilon}
\newcommand{\To}{\longrightarrow}
\newcommand{\C}{\mathcal{C}}
\newcommand{\K}{\mathcal{K}}
\newcommand{\T}{\mathcal{T}}
\newcommand{\bq}{\begin{equation}}
\newcommand{\eq}{\end{equation}}
\long\def\comments#1{ #1}
\comments{ }

\title[Conforming symmetric finite elements]{
A family of conforming mixed finite
   elements for linear elasticity on triangular grids
    }
\date{}
\author {Jun Hu}
\address{LMAM and School of Mathematical Sciences, Peking University,
  Beijing 100871, P. R. China.  hujun@math.pku.edu.cn}

\author {Shangyou Zhang}
\address{Department of Mathematical Sciences, University of Delaware,
    Newark, DE 19716, USA.  szhang@udel.edu }

\thanks{The first author was supported by  the NSFC Project 11271035,
    and  in part by the NSFC Key Project 11031006.}

\begin{abstract}
This paper presents a family of mixed finite elements on triangular grids for solving
the classical Hellinger-Reissner mixed problem of the elasticity equations.
In these elements,  the matrix-valued stress field is approximated by
 the full  $C^0$-$P_k$ space enriched by $(k-1)$ $H(\d)$  edge
  bubble functions on each internal edge, while
  the displacement field by the
    full discontinuous $P_{k-1}$ vector-valued space, for
  the polynomial degree $k\ge 3$.
  As a result,  compared with most of mixed   elements for linear elasticity
   in the literature, the basis of the stress space is surprisingly
    easy to construct.  The main challenge is to find the
    correct stress finite element space matching the full $C^{-1}$-$P_{k-1}$
    displacement space.
  The discrete stability analysis for the inf-sup condition
    does not rely on the usual Fortin  operator, which is difficult
    to construct.
  It is done by characterizing the divergence of local stress space which
    covers the $P_{k-1}$ space of displacement orthogonal to the local
    rigid-motion.
The  well-posedness condition and the optimal a priori error
estimate are proved for  this family of finite elements. Numerical
tests are presented to confirm the theoretical results.

  \vskip 15pt


 \vskip 15pt


\end{abstract}
\maketitle

\section{Introduction}
It is a challenge to design stable discretizations for the linear elasticity equations
   based on the
  Hellinger-Reissner variational principle, in which the stress and
  displacement are solved simultaneously.
This reason lies in, besides the usual  discrete K-ellipticity and
   B-B conditions, there is an additional symmetry constraint
      on the stress tensor for the problem under consideration.
   Many methods have been proposed to overcome this  difficulty, cf.  \cite{Amara-Thomas,
   Arnold-Brezzi-Douglas, Arnold-Douglas-Gupta,  Johnson-Mercier,
   Morley, Stenberg-1, Stenberg-2, Stenberg-3} for earlier works.
In a recent work  \cite{Arnold-Winther-conforming},  Arnold and Winther
   designed the first family of
    mixed finite element methods based on  polynomial shape  function spaces,
   which  was  motivated by a key
     observation: a discrete exact sequence
   guarantees the stability of the mixed method. From then on,
      various  stable  mixed elements have been constructed,
        see \cite{Adams-Cockburn,Arnold-Awanou,Arnold-Awanou-Winther,
    Arnold-Winther-conforming,Awanou, Chen-Wang},
     \cite{Arnold-Winther-n, Gopalakrishnan-Guzman-n,HuShi,
    Man-Hu-Shi, Yi-3D, Yi},  and   \cite{Arnold-Falk-Winther,
    Boffi-Brezzi-Fortin, Cockburn-Gopalakrishnan-Guzman,
    Gopalakrishnan-Guzman, Guzman}.
 Since most of these  elements  require a local commuting property
      which  implies that the usual Fortin operator
    can be constructed elementwise,   they  have many degrees
       of  freedom on each element such that they
     are not easy to be implemented; while the
    numerical examples can
    only be found in \cite{Carstensen-Eigel-Gedicke,
    Carstensen-Gunther-Reininghaus-Thiele2008} so far.

  In a recent paper,   a family of conforming
    mixed finite elements is proposed
    on rectangular grids
  for both  two and three dimensions.
  As a result the lowest order elements have $8$ plus $2$ and $18$
   plus $3$ degrees of freedom on each element for two and three dimensions,
    respectively, which are  simplifications of  two and three dimensional elements due to
  \cite{Hu-Man-Zhang2014}.
These elements were  motivated by an observation that  conformity of
the discrete  methods on  rectangular meshes  can be guaranteed by
 $H({\rm div})$-conformity of the  normal stress and
$H^1$-conformity of two corresponding variables for each  component
of the shear stress.
Such an idea was first explored in \cite{Hu-Man-Zhang2013} to  design the minimal mixed
 finite elements on rectangular grids in any dimension.
A new explicit constructional proof based on a macro-element technique
  was proposed to  show the discrete inf-sup condition for them.
In other words,  that  constructive proof avoids the discrete exact sequence
   of \cite{Arnold-Winther-conforming}, which is not possible  therein but used
   nearly everywhere \cite{Adams-Cockburn,Arnold-Awanou,Arnold-Awanou-Winther,
    Arnold-Winther-conforming,Awanou}.

 This paper presents a family of mixed finite elements on triangular grids.
In these elements,  the matrix-valued stress field is approximated by
 the full  $C^0$-$P_k$ space enriched by $(k-1)$ $H(\d)$
     edge bubble functions  on each internal edge, while
  the displacement field by a discontinuous vector-valued $P_{k-1}$ element for $k\ge 3$.
The main difficulty
  for the discrete stability analysis comes from the discrete inf-sup condition
  since it is impossible to  construct locally the usual Fortin  operator
    (for all $k\ge 3$).  To overcome  such a difficulty,   a new way of
  proof  is particularly   proposed
  to overcome it,  characterizing the divergence of local stress space which
    covers the $P_{k-1}$ space of displacement orthogonal to the local
    rigid-motion.
.

The new family of mixed elements is a simplification of the very first
   constructed family of stable elements of Arnold-Winther
   \cite{Arnold-Winther-conforming}.
For the $C^{-1}$-$P_{k-1}$ displace field,  the stress space of Arnold-Winther is
   the symmetric $H(\d)$-$P_{k+1}$ tensors whose divergence is in $P_{k-1}$,
    while ours is a subspace of symmetric $H(\d)$-$P_{k}$ tensors.
That is,  it is not needed to add those $P_{k+1}$ bubbles
  (of no approximation power) to the stress space,
   for the purpose of stability,  when $k\ge 3$.
 Computationally,  the new element is much simpler as there
   is no constraints on the polynomial degree deduction of divergence.
Mathematically, the new family of mixed element is the simplest one to achieve
  $P_{k-1}$ approximation for the displacement and
   $P_k$ approximation for the stress.
That is,  we eliminate all divergence-free stresses
  of no approximation power in the Arnold-Winther space.
As a result, the basis of our stress spaces is very  easy to construct. In fact,
 its basis can be directly derived by using the basis of the Lagrange element of order $k$.
However,  we do not improve the lowest order element in the Arnold-Winther family,
$k=2$, which will  be done, in a unified way,  with lower order elements for any space dimension,
in a forthcoming paper. We refer interested readers to \cite{Hu-Zhang2014b} for the extension
 to  the 3D case.

The rest of the paper is organized as follows. In Section 2, we
define the weak problem and the finite element method.
 In section 3,  we prove the well-posedness of the finite
    element problem, i.e. the discrete coerciveness and the
    discrete inf-sup condition.
  By which,  the optimal order convergence of the
   new element follows.
In Section 4, we provide some numerical results,
     using $P_3$, $P_4$ and $P_5$ finite elements and Arnold-Winther's $P_3$ element.

\section{The family of finite elements}

Based on the Hellinger-Reissner principle, the  linear elasticity
     problem within a stress-displacement ($\sigma$-$u$) form reads:
Find $(\sigma,u)\in\Sigma\times V :=H({\rm div},\Omega,\mathbb {S})
        \times L^2(\Omega,\mathbb{R}^2)$, such that
\an{\left\{ \ad{
  (A\sigma,\tau)+({\rm div}\tau,u)&= 0 && \hbox{for all \ } \tau\in\Sigma,\\
   ({\rm div}\sigma,v)&= (f,v) &\qquad& \hbox{for all \ } v\in V. }
   \right.\lab{eqn1}
}
Here the symmetric tensor space for stress $\Sigma$ and the
   space for vector displacement $V$ are, respectively,
  \an{   \lab{S}
  H({\rm div},\Omega,\mathbb {S})
    &:= \Big\{ \p{\sigma_{11} & \sigma_{12} \\\sigma_{21} &\sigma_{22}  }
     \in H(\d, \Omega)
    \ \Big| \ \sigma_{12} = \sigma_{21}  \Big\}, \\
     \lab{V}
     L^2(\Omega,\mathbb{R}^2) &:=
     \Big\{ \p{u_1 &  u_2}^T
          \ \Big| \ u_i \in L^2(\Omega) \Big\}  .}
	This paper denotes by $H^k(T,X)$ the Sobolev space consisting of
functions with domain $T\subset\mathbb{R}^2$, taking values in the
finite-dimensional vector space $X$, and with all derivatives of
order at most $k$ square-integrable. For our purposes, the range
space $X$ will be either $\mathbb{S},$ $\mathbb{R}^2,$ or
$\mathbb{R}$.
$\|\cdot\|_{k,T}$ is the norm of $H^k(T)$. $\mathbb{S}$ denotes
the space of symmetric tensors, $H({\rm div},T,\mathbb{S})$
consists of square-integrable symmetric matrix fields with
square-integrable divergence. The H(div) norm is defined by
$$\|\tau\|_{H({\rm div},T)}^2:=\|\tau\|_{L^2(T)}^2
   +\|{\rm div}\tau\|_{L^2(T)}^2.$$
$L^2(T,\mathbb{R}^2)$ is the space of vector-valued functions
which are square-integrable.

Throughout the paper, the compliance tensor
$A=A(x):\mathbb{S}~\rightarrow~\mathbb{S}$, characterizing the
properties of the material, is bounded and symmetric positive
definite uniformly for $x\in\Omega$.

This paper deals with a pure displacement problem \meqref{eqn1} with the
  homogeneous boundary condition that $u\equiv 0$ on
  $\partial\Omega$.
But the method and the analysis work for mixed boundary value problems
   and the pure traction boundary problem.

The  domain $\Omega$ is subdivided by a family of quasi-uniform
  triangular grids  $\mathcal{T}_h$ (with the grid size $h$).
We introduce the finite element space of order $k$ ($k\ge3$)  on $
   \mathcal{T}_h$.
The displacement space is the full $C^{-1}$-$P_{k-1}$ space
 \an{ \lab{Vh}
   V_h = \{v\in L^2(\Omega,\mathbb{R}^2),
         v|_K\in P_{k-1}(K, \mathbb{R}^2)\ \hbox{ for all } K\in\mathcal{T}_h \}.
    }

The stress space is the full $C^0$-$P_k$ space enriched by $(k-1)$ $H(\d)$
   edge bubble functions on each internal edge.
We define the edge bubble functions first.
Let $\triangle \b x_0\b x_1 \b x_2=:K\in \mathcal{T}_h$ with three edges $E_i$ and corresponding three
   barycentric variables $\lambda_i$.
Here $\lambda_i$ is a linear function which vanishes on edge $E_i$ and
   assumes a nodal value 1 at the opposite vertex $\b x_i$,  see Figure
   \mref{bary}. Given $E_i=\vec{\b x_{i-1}\b x_{i+1}}$, its two endpoints are  $\b x_{i-1}$ and $\b x_{i+1}$,
    which allows for defining   its $k-1$ interior nodal points by
  \begin{equation}\label{E}
  \b x_{E_i,j}=\frac{j}{k}\b x_{i-1}+\frac{k-j}{k}\b x_{i+1}, j=1, \cdots, k-1.
  \end{equation}
  We also define $\frac{(k-1)(k-2)}{2}$   nodal points inside $K$ by
   \begin{equation}\label{K}
   \b x_{K, l, m}=\frac{l}{k}\b x_0+\frac{m}{k}\b x_1+\frac{k-l-m}{k}\b x_2, 1\leq l,  m \text{ and }l+m\leq k-1.
   \end{equation}
    Then the nodes for the Lagrange element of order $k$ is
    \begin{equation*}
    \begin{split}
    \b \mathbb{X}_K=\{\b x_i, i=0,1,2\}&\cup\{\b x_{E_i,j}, i=0,1,2, j=1, \cdots, k-1\}\\
    &\cup \{\b x_{K,l, m}, 1\leq l, m \text{ and } l+m\leq k-1\}.
  \end{split}
    \end{equation*}
    Given node $\b x_{E_i,j}$ on edge $E_i$, $j=1, \cdots, k-1$,
     let $\phi_{E_i, j}\in P_k(K, \mathbb{R})$
    be its associated  nodal basis  function of the Lagrange element of order $k$ such that
     \begin{equation}\label{basis}
     \phi_{E_i, j}(\b x_{E_i, j})=1 \text{ and } \phi_{E_i, j}(\b x^\prime)=0
     \text{ for any }\b x^\prime\in \b \mathbb{X}_K \text{ other than }\b x_{E_i, j}.
     \end{equation}
 Let $\b n_i=\langle n_{i,1}, n_{i,2}\rangle$ and $\b n_i^\perp=\langle -n_{i,2}, n_{i,1}\rangle$ be normal and tangent vectors on
   edge $E_i$,  respectively,  see Figure \mref{affine}. We define a matrix of rank one by
   \begin{equation}\label{rankone}
   \mathbb{T}_{E_i}=\b n_i^\perp {\b n_i^\perp}^T.
   \end{equation}
   With these $(k-1)$ edge bubble functions $\phi_{E_i, j}$ on each edge and the matrix $\mathbb{T}_{E_i}$ of rank one,  we can define   exactly $(k-1)$ stress functions $\tau_{\mathbb{E}_i,j}$ by
 \an{\lab{z-f}
  \tau_{E_i,j}=\phi_{E_i,j}\mathbb{T}_{E_i} , \quad j=1,2,\dots,k-1, \ i=0,1,2.}

  By the  definition, we have
  \begin{equation}
  \tau_{E_i,j}\cdot\b n_l|_{E_l}=0, i,l=0,1,2, j=1, \cdots, k-1,
  \end{equation}
  which implies that they are $H(\d )$ bubble functions on element $K$.  We define
  \begin{equation}
  \Sigma_{\partial K, b}=\text{span}\{\tau_{E_i, j}, i=0, 1, 2, j=1, \cdots, k-1\}.
  \end{equation}

  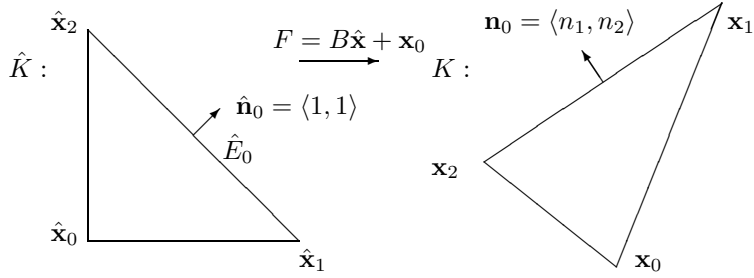
\begin{figure}[htb]\setlength\unitlength{1pt}\begin{center}
    \begin{picture}(300,110)(0,0)
  \put(20,10){\begin{picture}(120,80)(0,0)\put(-20,72){$\hat K:$}
      \put(10,10){\line(1,0){80}}  \put(10,10){\line(0,1){80}}
      \put(10,90){\line(1,-1){80}}
       \put(-4,10){$\hat{\b x}_0$}  \put(90,0){$\hat{\b x}_1$}
     \put(-4,90){$\hat{\b x}_2$} \put(61,41){$\hat E_0$}
      \put(50,50){\vector(1,1){10}} \put(66,58){$\hat{\b n}_0=\langle 1, 1 \rangle$}
    \end{picture}  }
   \put(110,88){\vector(1,0){30}}   \put(100,93){$F=B\hat{\b x}+\b x_0$}

  \put(180,10){\begin{picture}(120,80)(0,0) \put(-20,72){$K:$}
      \put(0,40){\line(5,-4){50}}  \put(50,0){\line(2,5){40}}
      \put(90,100){\line(-3,-2){90}}  \put(45,70){\vector(-2,3){8}}
      \put(-20,35){$\b x_2$}  \put(57,0){$\b x_0$}
      \put(92,90){$\b x_1$}   \put(0,90){$\b n_0=\langle n_1, n_2 \rangle$}
    \end{picture}  }

   \end{picture}
\caption{\lab{affine} A reference triangle and a general triangle with an edge
    normal vector.  }
    \end{center} \end{figure}

 \begin{figure}[htb]\setlength\unitlength{0.95pt}\begin{center}
    \begin{picture}(330,100)(0,0)
  \put(10,10){\line(1,0){80}}  \put(10,10){\line(0,1){80}}
  \put(10,90){\line(1,-1){80}}
   \put(-4,10){$\b x_0$} \put(-14,0){$\lambda_0(\b v_0)=1$} \put(90,0){$\b x_1$}
 \put(-4,90){$\b x_2$} \put(51,51){$E_0$}\put(31,71){$\lambda_0(E_0)=0$}
   \put(-5,50){$E_1$} \put(40,0){$E_2$}
  \def\trs{  \put(10,10){\line(1,0){80}}  \put(10,10){\line(0,1){80}}
  \put(10,90){\line(1,-1){80}}
   \multiput(10,10)(26.6667,0){4}{\circle{3}}
   \multiput(10,36.667)(26.6667,0){3}{\circle{3}}
   \multiput(10,63.3333)(26.6667,0){2}{\circle{3}}
    \put(10,90){\circle{3}}  }
   \put(100,0){\begin{picture}(100,100)\trs  \put(36.666,10){\line(0,1){53.3333}}
    \put(63.333,36.6667){\circle*{3}} \put(63,40){$\phi_{0,1}=1$}
    \put(20,82){$\tilde \phi_{0,1}=\lambda_1(\lambda_1-\frac 13)\lambda_2$}
    \end{picture}  }
  \put(230,0){\begin{picture}(100,100)\trs
     \put(10,36.666){\line(1,0){53.3333}}
    \put(36.6667,63.3333){\circle*{3}}\put(40,63){$\phi_{0,2}=1$}
    \put(20,82){$\tilde \phi_{0,2}=\lambda_1(\lambda_2-\frac 13)\lambda_2$}
    \end{picture}  }
   \end{picture}
\caption{\lab{bary}The barycentric variables $\lambda_i$ (linear functions)
        on $K=\b x_0\b x_1\b x_2$,
      and two edge bubble functions of $P_3$ on edge $E_0$, where
        $\phi_{0,1}=\tilde\phi_{0,1}/\tilde\phi_{0,1}
	   (\lambda_1=\frac 23,\lambda_2=\frac 13)$. }
    \end{center} \end{figure}
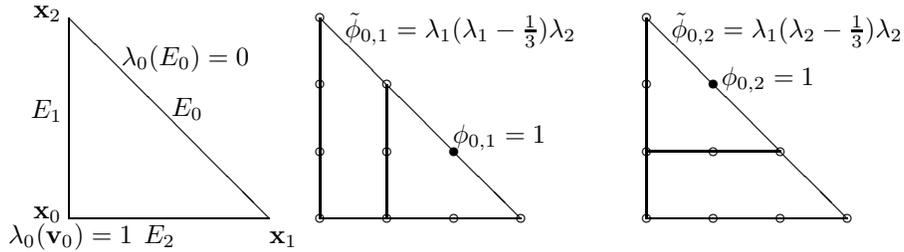

The finite element space of order $k$ ($k\ge 3$)
     for the stress approximation is
\an{\lab{Sh}
 \Sigma_h =\Big\{~\sigma
    &\in H(\d,\Omega,\mathbb{S}),  \sigma=\sigma_c+\sigma_b,
             \ \sigma_c\in H^1(\Omega, \mathbb{S}), \\
     & \ \sigma_c|_K\in P_k(K, \mathbb{S})
	        \,,
	\nonumber  \ \sigma_b|_K\in\Sigma_{\partial K, b},  \forall K\in\mathcal{T}_h \Big\}, }
 which is a $H(\d )$ bubble enrichment of the $H^1$ space
  \an{\lab{wSh}
 \widetilde \Sigma_h =\Big\{~\sigma
    \in H(\d,\Omega,\mathbb{S}),
           & \ \sigma \in H^1(\Omega,\mathbb{S}), \ \sigma|_K\in P_k(K, \mathbb{S})
	        \ \forall K\in\mathcal{T}_h \Big\}.
            }
To define a basis of $\Sigma_h$, we need the orthogonal complement matrices $\mathbb{T}_{E, j}^\perp\in\mathbb{S}$, $j=1, 2$, of  matrix $\mathbb{T}_E$ for any edge $E$ of $\mathcal{T}_h$,  which are defined by
\begin{equation}
\mathbb{T}_{E, j}^\perp:\mathbb{T}_E=0, \mathbb{T}_{E, j}^\perp: \mathbb{T}_{E, j}^\perp=1,
 \text{ and }\mathbb{T}_{E, 1}^\perp: \mathbb{T}_{E, 2}^\perp=0,
\end{equation}
where the inner product $A:B=a_{11}b_{11}+a_{12}b_{12}+a_{21}b_{21}+a_{22}b_{22}$ for two
 matrices $A=\{a_{ij}\}_{i,j=1}^2$ and $B=\{b_{ij}\}_{i,j=1}^2$. The canonical basis of $\mathbb{S}$
  reads
  \begin{equation}
  \mathbb{T}_1=\p{1 & 0\\ 0&0}, \mathbb{T}_2=\p{0& 1\\1 &0}, \text{ and }\mathbb{T}_3=\p{0&0\\0 &1}.
  \end{equation}
  Let $\mathcal{X}_{\mathbb{E}}$ denote  all interior  nodes, defined in \eqref{E},  of all the  edges,  $\mathcal{X}_{\mathbb{K}}$ denote all interior nodes, defined in \eqref{K}, of all the elements,  and   $\mathcal{X}_{\mathbb{V}}$ denote all the  vertices of  $\mathcal{T}_h$.  Define the Lagrange element
   space of order $k$ by
   $$
   \mathbb{P}_h:=H^1(\Omega, \mathbb{R})\cap \{v\in L^2(\Omega), v|_K\in P_k(K, \mathbb{R}), \forall K\in\mathcal{T}_h\}.
   $$
   Given node $\b x\in \mathcal{X}_{\mathbb{V}}\cup \mathcal{X}_{\mathbb{E}}\cup \mathcal{X}_{\mathbb{K}}$,
    let $\phi_{\b x}\in \mathbb{P}_h$  be its associated  nodal basis function, which is similarly defined
     as $\phi_{E_i, j}$ in \eqref{basis}.
  The basis functions of $\Sigma_h$ can be classified into  four classes:
 \begin{enumerate}
 \item Vertex--based basis functions:  given vertex $\b x\in \mathcal{X}_{\mathbb{V}}$,  its three associated basis functions of  $\Sigma_h$ read
     $$
     \tau_{V, \b x, i}=\phi_{\b x}\mathbb{T}_i, i=1, 2, 3.
     $$
 \item Volume--based basis functions: given node $\b x \in \mathcal{X}_{\mathbb{K}}$ inside $K$,  its three associated basis functions of  $\Sigma_h$ read
     $$
     \tau_{K, \b x, i}=\phi_{\b x}\mathbb{T}_i, i=1, 2, 3.
     $$
   \item Edge--based basis functions with nonzero fluxes: given node $\b x \in \mathcal{X}_{\mathbb{E}}$ on edge $E$, its two  associated basis functions with nonzero fluxes of  $\Sigma_h$ read
 $$
     \tau_{E, \b x, i}^{(nb)}=\phi_{\b x}\mathbb{T}_{E, i}^\perp, i=1, 2.
     $$
 \item Edge--based bubble functions: given node $\b x \in \mathcal{X}_{\mathbb{E}}$ on edge $E$ which is
  shared by elements $K_1$ and $K_2$, its bubble functions in  $\Sigma_h$ read
   $$
     \tau_{E, \b x, i}^{(b)}=\phi_{\b x}|_{K_i}\mathbb{T}_{E}, i=1, 2.
    $$

 \end{enumerate}
It is  straightforward to see that these functions defined in the above four terms form a basis
 of $\Sigma_h$, which are   very easy to  construct.

The mixed finite element approximation of Problem (1.1) reads: Find
   $(\sigma_h,~u_h)\in\Sigma_h\times V_h$ such that
 \e{ \left\{ \ad{
    (A\sigma_h, \tau)+({\rm div}\tau, u_h)&= 0 &&
              \hbox{for all \ } \tau \in\Sigma_h,\\
     (\d\sigma_h, v)& = (f, v) &&  \hbox{for all \ } v\in V_h.
      } \right. \lab{DP}
    }
 It follows from the definition of $V_h$ ($P_{k-1}$ polynomials)
    and $\Sigma_h$ ($P_k$ polynomials) that
   \a{ \d  \Sigma_h \subset V_h.}
This, in turn, leads to a strong divergence-free space:
 \an{ \lab {kernel}
    Z_h&:= \{\tau_h\in\Sigma_h \ | \ (\d\tau_h, v)=0 \quad
	\hbox{for all } v\in V_h\}\\
    \nonumber
          &= \{\tau_h \in\Sigma_h \ | \  \d \tau_h=0
    \hbox{\ pointwise } \}.
    }

\section{Stability and convergence}
The convergence of the finite element solutions follows
   the stability and the standard approximation property.
So we consider first the well-posedness  of the discrete problem
    \meqref{DP}.
By the standard theory,  we only need to prove
   the following two conditions, based on their counterpart at
    the continuous level.

\begin{enumerate}
\item K-ellipticity. There exists a constant $C>0$, independent of the
   meshsize $h$ such that
    \an{ \lab{below} (A\tau, \tau)\geq C\|\tau\|_{H(\d)}^2\quad
       \hbox{for all } \tau \in Z_h, }
    where $Z_h$ is the divergence-free space defined in \meqref{kernel}.

\item  Discrete B-B condition.
    There exists a positive constant $C>0$
            independent of the meshsize $h$, such that
    \an{\lab{inf-sup}
   \inf_{0\neq v\in V_h}   \sup_{0\neq\tau\in\Sigma_h}\frac{({\rm
        div}\tau, v)}{\|\tau\|_{H(\d)}  \|v\|_{L^2(\Omega)} }\geq
    C .}
\end{enumerate}

It follows from $\d  \Sigma_h \subset V_h$ that $\d  \tau=0$ for
   any $\tau\in Z_h$. This implies the above K-ellipticity condition
	\meqref{below}.

It remains to show the discrete B-B condition \meqref{inf-sup},
  in the following two lemmas.

\begin{lemma}\label{lemma1}
For any $v_h\in V_h$,  there is a $\tau_h \in \widetilde \Sigma_h$  such that,
  for any polynomial $p\in P_{k-2}(K, \mathbb{R}^2)$,
   \bq\label{l-1}  \int_K (\d\tau_h-v_h) \cdot p\,d\b x=0
      \quad \hbox{\rm
      and } \quad \|\tau_h\|_{H(\d)}\leq C\|v_h\|_{L^2(\Omega)}. \eq
     \end{lemma}

\begin{proof}
  Let $v_h\in V_h$.
  By the stability of the continuous formulation,
      cf. \cite{Arnold-Winther-conforming}, there is
    a $\tau \in \Sigma\cap H^1(\Omega, \mathbb{S})$ such that,
   \a{ \d\tau=v_h \quad \hbox{\rm
      and } \quad \|\tau\|_{H^1(\Omega)}\le  C\|v_h\|_{L^2(\Omega)}. }
As $\tau \in H^1(\Omega, \mathbb{S})$,  we modify
   the Scott-Zhang \cite{Scott-Zhang}
   interpolation operator slightly to define a flux preserving
   interpolation.
 \a{  I_h \ : \ \Sigma\cap H^1(\Omega, \mathbb{S})
              & \to \Sigma_h \cap H^1(\Omega, \mathbb{S}) =\widetilde \Sigma_h , \\
         \ \tau=\begin{pmatrix}\tau_{11} & \tau_{12}\\ \tau_{12} & \tau_{22} \end{pmatrix}
      & \mapsto \tau_h = \begin{pmatrix}\tau_{11,h} & \tau_{12,h}\\
      \tau_{12,h} & \tau_{22,h} \end{pmatrix}= I_h \tau. }
Here the interpolation is done inside a subspace,
   the continuous finite element subspace
     $\Sigma_h \cap H^1(\Omega, \mathbb{S})$.
$I_h \tau$ is defined by its values at the Lagrange nodes.

At a vertex node $\b x_i$,
   $I_h\tau (\b x_i)$ is defined as the nodal value of $\tau$ at the
     vertex if $\tau$ is continuous, but in general, $I_h\tau(\b x_i)$ is
       defined as an average value on an edge at the vertex,
        as in \cite{Scott-Zhang}.
After defining the nodal values at vertices of triangles,
  the nodal values of $\tau_h$
      at the nodes inside each edge are defined by
   the $L^2$-orthogonal projection on the edge:
  \an{ \lab{a-1}
   \forall p\in P_{k-2}(E, \mathbb{R}), \ \left\{ \ad{
    \int_E \tau_{h,11}  p  \, ds &=
    \int_E \tau_{ 11}  p  \, ds, \\
    \int_E  \tau_{h,12}  p  \, ds &=
    \int_E  \tau_{ 12}  p  \, ds, \\
     \int_E  \tau_{h,22}  p  \, ds &=
    \int_E  \tau_{ 22}  p  \, ds, }\right. }
 where $E$ is an edge in the triangulation $\mathcal{T}_h$.
At the Lagrange nodes inside triangles,
   $I_h\tau$ is defined by
   the $L^2$-orthogonal projection on the triangle:
 \an{ \lab{a-2}
    \int_K  \tau_{ij,h}  p  \, d\b x &= \int_K \tau_{ij}  p d\b x \quad
   \forall p\in P_{k-3}(K, \mathbb{R}),   }
    where $K$ is an element of  $\mathcal{T}_h$.
  It follows by the stability of the Scott-Zhang operator that
   \a{ \| I_h \tau \|_{H(\d)}\le C \| \tau \|_{H^1(\Omega)}\le C\|v_h\|_{L^2(\Omega)}. }
By \meqref{a-1} and \meqref{a-2},
  we get the a partial-divergence matching property of $I_h$:
     for any $p\in P_{k-2}(K, \mathbb{R}^2)$,
 \a{  \int_K (\d\tau_h-v_h)\cdot p\, d\b x
    & = \int_{\partial K}  \tau_h \b n\cdot p \,ds
        -\int_K \tau_h: \nabla p\, d \b x
           - \int_K v_h\cdot p\,d\b x \\
        & = \int_{\partial K}  \tau \b n \cdot p\, ds
        -\int_K \tau: \nabla p\, d \b x
           - \int_K v_h\cdot p \,d\b x \\
	 	 &= \int_K (\d\tau -v_h)\cdot p\, d\b x = 0. }
\end{proof}

\begin{lemma}\lab{lemma2}
For any $v_h\in V_h$, if
  \an{\lab{2o}
   \int_K  v_h \cdot p d\b x=0 \quad\hbox{ for all $p\in P_{k-2}(K, \mathbb{R}^2)$, } }
   there is a $\tau_h \in \Sigma_h$  such that
   \an{ \lab{l-2}  \d\tau_h = v_h
      \quad \hbox{\rm
      and } \quad \|\tau_h \|_{H(\d)}\leq C\|v_h\|_{L^2(\Omega)}. }
     \end{lemma}
\begin{proof} We first define the local spaces of bubble stress functions.
Let $\b x_0$, $\b x_1$ and $\b x_2$ be the three vertices of a triangle $K$.
The referencing mapping is then, cf. Figure \mref{affine},
\a{  \b x &= F_K(\hat {\b x})
     = \b x_0 + \p{ \b x_1 -\b x_0 &
                     \b x_2 -\b x_0 } \hat{\b x}. }
   Then
\an{ \lab{imap}  \hat {\b x} &= \p{\b n_1^T\\ \b n_2^T } (\b x-\b x_0), }
where \a{ \p{\b n_1^T\\ \b n_2^T }=\p{ \b x_1 -\b x_0 &
                     \b x_2 -\b x_0 }^{-1}. }
Due to the inverse matrix relation,  these two vectors $\b n_1, \b n_2$ are
  orthogonal to edges   $\vec{\b x_0\b x_2}$ and $\vec{\b x_0\b x_1}$, respectively.
By \meqref{imap}, they are coefficients of the barycentric variables:
  \a{ \lambda_1 &=  \b n_1 \cdot (\b x-\b x_0), \\
     \lambda_2 &= \b n_2 \cdot (\b x-\b x_0), \\
     \lambda_0 &= 1-\lambda_1-\lambda_2. }
With them,  we define the $H(\d, K, {\mathbb S})$ bubble functions
\an{\lab{bh}  \Sigma_{K,b} &= \operatorname{span}
    \{ \lambda_2\lambda_0p_1
                   \b n_1^\perp {\b n_1^\perp}^T,
       \lambda_0\lambda_1 p_2  \b n_2^\perp {\b n_2^\perp}^T,
     \lambda_1\lambda_2p_0  \b n_0^\perp {\b n_0^\perp}^T \},
              }
        where $p_{1}, p_{2}$ and
          $p_0\in P_{k-2}(K, \mathbb{R})$, and
      \a{ \b n_1^\perp &=\p{-n_{12} \\ n_{11}}, \quad\hbox{ if \ }
                     \b n_1  = \p{ n_{11} \\ n_{12}}, \\
          \b n_0 & = -\b n_1 -\b n_2. }
  Note that $ \tau_h \cdot \b n_j =\b 0$  on the three edges ($\lambda_j=0$),
    for all $\tau_h \in  \Sigma_{K,b}$.
 Thus, the match of $\d \tau_h = v_h$ is done locally on $K$, independently
   of the matching on next element.

We begin to prove the lemma.
Let $v_h\in V_h$ satisfying \meqref{2o}.
We show there is a local $\tau \in  \Sigma_{K,b} $ such that $\d \tau  = v_h$,
  on each element $K$.
As $v_h$ satisfies \meqref{2o}, $v_h|_K \in  V_{K,\perp}$ where
   $V_{K,\perp}$ is the rigid-motion free space
   \a{ V_{K,\perp} = \{ v_h \in P_{k-1}(K, \mathbb{R}^2),  \int_K v_h \cdot \p{a-by \\ c+bx} d\b x=\b 0,
          \ \forall a,b,c \in {\mathbb R} \}. }
We prove $\d \Sigma_{K,b} = V_{K,\perp}$ next.
By definition,  $\d \Sigma_{K,b} \subset V_{K,\perp}$.
If $\d \Sigma_{K,b} \ne   V_{K,\perp}$, there is a $v_h\in V_{K,\perp}$
   orthogonal to $\d \Sigma_{K,b}$, i.e.,
  \a{
   \int_K  \d \tau \cdot v_h d\b x =- \int_K \tau : \epsilon(v_h) d\b x=0
    \quad \forall \tau \in \Sigma_{K,b}, }
  where $\epsilon(v_h)$ is the symmetric gradient, $(\nabla v_h + \nabla^T v_h)/2$.
We show next $v_h=0$.

Let $\{\mathbb{M}_i, i=0,1,2\}$ be the dual basis,  of
  \a{ \mathbb{T}_i= \b n_i^\perp {\b n_i^\perp}^T, \quad i=1,2,0, }
i.e.
  \a{ \mathbb{M}_j: \mathbb{T}_i
    =\delta_{ij}. }
As noted above,  $\{\b n_i\}$ are three normal vectors of element $K$.
In Lemma \ref{lemma3} below, we shall prove that the three matrices $\mathbb{T}_i$, $i=0, 1, 2$, are linearly independent.
Hence the above equation has a unique solution $\mathbb{M}_j$. Therefore we have a unique expansion, as $\epsilon(v_h) \in P_{k-2}(K,\mathbb{S})$,
  \a{ \epsilon(v_h)= q_1 \mathbb{M}_1 + q_2 \mathbb{M}_2 + q_0 \mathbb{M}_0, \quad \hbox{
      for some } \ q_1,q_2,q_0 \in P_{k-2}(K, \mathbb{R}).  }
Selecting $\tau_1 = \lambda_2\lambda_0q_1
                   \b n_1^\perp {\b n_1^\perp}^T\in \Sigma_{K,b}$,
         we have
   \a{ 0 = \int_K \tau_1 : \epsilon(v_h) d\b x
         = \int_K \lambda_2\lambda_0 q_1^2(\b x) d\b x. }
As $\lambda_2\lambda_0>0$ on $K$,  we conclude that $q_1=0$.
Similarly, $q_2$ and $q_0$ are zero.  Since  $v_h=0$ is not a rigid motion, this shows that $v_h=0$.

As we assume $k\ge 3$,  the condition \meqref{2o}
     that $v_h$ is orthogonal to $P_{k-2}(K, \mathbb{R}^2)$
   implies that $v_h\in V_{h,\perp}=\d \Sigma_{K,b}.$
The selection of $\tau_h$, locally on element $K$,  is made by
   \a{ \|\tau_h\|_{L^2(\Omega)} =\min \{ \|\tau\|_{L^2(\Omega)},
     \d \tau =v_h, \ \tau\in \Sigma_{K,b}\}. }
The boundedness of $\d$ operator in \meqref{l-2} follows the scaling argument
   with affine mappings.
Thus \meqref{l-2} holds.
\end{proof}

\begin{lemma}\label{lemma3} Let $\b v_1\in \mathbb{R}^2$ and $\b v_2\in \mathbb{R}^2$, and $\b v_0=\b v_1+\b v_2$. Suppose that
$\b v_1$ and $\b v_2$ are linearly  independent.  Then three matrices $\b v_1 \b v_1^T$, $\b v_2 \b v_2^T$, and $\b v_0 \b v_0^T$ are linearly independent.
\end{lemma}
\begin{proof}  Let $\b v_1=(a_1, a_2)^T$ and $\b v_2=(b_1, b_2)^2$, which  leads to  $\b v_0=(a_1+b_1, a_2+b_2)^T$. Hence
$$
\b v_1 \b v_1^T=\begin{pmatrix}a_1^2 & a_1a_2\\ a_1a_2 & a_2^2 \end{pmatrix}, \b v_2 \b v_2^T=\begin{pmatrix}b_1^2 & b_1b_2\\ b_1b_2 & b_2^2 \end{pmatrix} $$
and
$$
\b v_0 \b v_0^T
=\begin{pmatrix}(a_1+b_1)^2 & (a_1+b_1)(a_2+b_2)\\ (a_1+b_1)(a_2+b_2) &(a_2+ b_2)^2 \end{pmatrix}
$$
To prove the desired result, it suffices to
 show that the rank of the  matrix
 $$
  \begin{pmatrix}a_1^2 & b_1^2 & (a_1+b_1)^2\\a_2^2 & b_2^2 & (a_2+b_2)^2 \\a_1a_2& b_1b_2&(a_1+b_1)(a_2+b_2)\end{pmatrix},
 $$
 is three.
 A direct calculation finds that the determinant of the above matrix is $(a_1b_2-a_2b_1)^3$. Since $\b v_1$ and $\b v_2$ are linearly independent, we have
 $$
 a_1b_2-a_2b_1\not=0,
 $$
which completes the proof.
\end{proof}

\begin{remark} The lemma \ref{lemma2} can be proved differently,  by counting the dimension of
   vector spaces.
 Due to the linearly independent vectors (in matrix form),
  \a{ \dim  \Sigma_{K,b}  &= 3 \dim P_{k-2} = \frac 32 k^2 - \frac 32 k.}
If we can show that the div-free bubbles of $\Sigma_{K,b}$ must be the
   bubble Airy functions, namely,
   \an{\lab{af} \d \tau_h =0 \ \Rightarrow \
        \tau_h = \p{ \frac{\partial^2  w_h}{\partial y^2} &
                    -\frac{\partial^2 w_h }{\partial x\partial y} \\
                    -\frac{\partial^2 w_h }{\partial x\partial y} &
                    \frac{\partial^2  w_h}{\partial x^2}}
    } where $w_h=w_K b_K^2$ for some $w_K\in P_{k-4}(K, \mathbb{R})$, where $b_K=\lambda_0\lambda_1\lambda_2$ is
   the element $P_3$ bubble,
then we would get
  \a{ \dim \d  \Sigma_{K,b}  &=  \dim  \Sigma_{K,b}- \dim  P_{k-4} \\
                     & =    k^2 + k - 3
                     = 2 \dim P_{k-1} -3 =\dim V_{K,\perp}. }
 As $\d \Sigma_{K,b}\subset V_{K,\perp}$, the dimension counting will prove
     $\d \Sigma_{K,b}=V_{K,\perp}$.

 We are going to prove \meqref{af}.
  Since $w_h$ can be selected up to a linear function,
  we start to take $w_h$ such that it vanishes at three vertices of element $K$.
   Since $\tau_h\in \Sigma_{K,b}$,  it follows that
 $$
 \frac{\partial }{\partial \b n_i^\perp}\frac{\partial  w_h}{\partial x}|_{E_i}
 =\frac{\partial }{\partial \b n_i^\perp}\frac{\partial  w_h}{\partial y}|_{E_i}
  =0, i=0, 1, 2.
 $$
This implies that
 \begin{equation}\label{tang}
 \frac{\partial ^2 w_h}{\partial (\b n_i^\perp)^2}|_{E_i}
 =\frac{\partial ^2 w_h}{\partial \b n_i^\perp \partial \b n_i}|_{E_i}=0, i=0, 1, 2.
 \end{equation}
Hence  $\frac{\partial w_h}{\partial \b n_i^\perp}$ is a constant on $E_i$.
  Since $w_h$ vanishes on three vertices of $K$,
   this indicates that $w_h$ vanishes on $E_i$,
   which implies that $\frac{\partial w_h}{\partial \b n_i^\perp}=0$ on $E_i$.
   Consequently, $\nabla w_h$ vanishes on three vertices of $K$.
   By \eqref{tang}, $\frac{\partial  w_h}{ \partial \b n_i}$ is a constant on $E_i$.
  This implies that $\frac{\partial  w_h}{ \partial \b n_i}=0$ on $E_i$,
   which  completes the proof of  \meqref{af}.
\end{remark}

We are in the position to show the well-posedness of the discrete problem.
\begin{lemma}
 For the discrete problem (\ref{DP}), the K-ellipticity \meqref{below}
    and the discrete B-B
 condition \meqref{inf-sup} hold uniformly.
  Consequently,  the discrete
     mixed problem \meqref{DP} has a unique solution
         $(\sigma_h,~u_h)\in\Sigma_h\times V_h$.
\end{lemma}
\begin{proof}  The  K-ellipticity immediately follows from the fact that
   $\d  \Sigma_h \subset V_h.$ Therefore we only need to
 prove the  discrete B-B  condition \meqref{inf-sup}. For any $v_h\in V_h$,
 it follows from Lemma \ref{lemma1} that there exists a $\tau_{1}\in \Sigma_h$ such that,
  for any polynomial $p\in P_{k-2}(K, \mathbb{R}^2)$,
   \bq  \int_K (\d\tau_1-v_h) \cdot p\,d\b x=0
      \quad \hbox{\rm
      and } \quad \|\tau_1\|_{H(\d)}\leq C\|v_h\|_{L^2(\Omega)}. \eq
Then it follows from Lemma \ref{lemma2} that
 there is a $\tau_2 \in \Sigma_h$  such that
   \an{   \d\tau_2 = v_h-\d\tau_1
      \quad \hbox{\rm
      and } \quad \|\tau_2 \|_{H(\d)}\leq C\|\d\tau_1-v_h\|_{L^2(\Omega)},}
Let $\tau=\tau_1+\tau_2$.  This implies that
\begin{equation}
\d\tau=v_h \text{ and } \|\tau\|_{H(\d)}\leq C\|v_h\|_{L^2(\Omega)},
\end{equation}
this proves the discrete B-B condition \meqref{inf-sup}.
\end{proof}

\begin{theorem}\label{MainError} Let
  $(\sigma, u)\in\Sigma\times V$ be the exact solution of
   problem \meqref{eqn1} and $(\tau_h, u_h)\in\Sigma_h\times
   V_h$ the finite element solution of \meqref{DP}.  Then,
   for $k\ge 3$,
\an{ \lab{t1} \|\sigma-\sigma_h\|_{H({\rm div})}
    + \|u-u_h\|_{L^2(\Omega)}&\le     Ch^k(\|\sigma\|_{H^{k+1}(\Omega)}+\|u\|_{H^k(\Omega)}).
      }
\end{theorem}
\begin{proof}
 The stability of the elements and the standard theory of mixed
  finite element methods \cite{Brezzi, Brezzi-Fortin} give the
  following quasi--optimal error estimate immediately
\an{
  \label{theorem-err1} \|\sigma-\sigma_h\|_{H({\rm
  div})}+\|u-u_h\|_{L^2(\Omega)}\leq C \inf\limits_{\tau_h\in\Sigma_h,v_h\in
  V_h}\left(\|\sigma-\tau_h\|_{H({\rm div})}+\|u-v_h\|_{L^2(\Omega)}\right).}
Let $P_h$ denote the local $L^2$ projection operator,
   or triangle-wise interpolation operator,  from $V$ to $V_h$,
  satisfying the error estimate
\an{\label{proj-error}
   \|v-P_hv\|_{L^2(\Omega)}\leq Ch^k\|v\|_{H^k(\Omega)} \text{ for any }v\in H^k(\Omega, \mathbb{R}^2). }
Choosing $\tau_h=I_h\sigma\in \Sigma_h$
    where $I_h$ is defined in \meqref{a-1} and \meqref{a-2},
we have \cite{Scott-Zhang},  as $I_h$ preserves symmetric $P_k$ functions locally,
   \an{ \lab{p-err2}
       \|\sigma -\tau_h\|_{L^2(\Omega)} + h |\sigma -\tau_h|_{\d}
        \le Ch^{k+1} \|\sigma\|_{H^{k+1}(\Omega)}. }
 Let $v_h= P_h v$ and $\tau_h=I_h\sigma$ in (\ref{theorem-err1}),
  by (\ref{proj-error}) and \meqref{p-err2}, we
   obtain  \meqref{t1}.
\end{proof}

\begin{remark}
To prove an optimal error estimate for the stress in the $L^2$ norm, we  can follow the idea from \cite{Stenberg-1}  to use
 a mesh dependent norm technique.  In particular, this will lead to
 $$
 \|\sigma-\sigma_h\|_{L^2(\Omega)}\leq Ch^{k+1}|\sigma|_{H^{k+1}(\Omega)}.
 $$
\end{remark}

\begin{remark}
The extension to nearly incompressible or incompressible elastic materials is possible. In the homogeneous
isotropic case the compliance tensor is given by
\a{
      A \tau &= \frac 1{2\mu} \left(
       \tau - \frac{\lambda}{2\mu + n \lambda} \operatorname{tr}(\tau)
        \delta \right), \quad n=2,  }
  where $\delta=\p{1 &0\\0&1}$, and  $\mu > 0$, $\lambda>0$ are the Lam\'{e} constants. For our mixed
method, as for most methods based on the Hellinger--Reissner principle,
one can prove that the error estimates hold uniformly in $\lambda$.
In the analysis above we use the fact that
$$
\alpha\|\tau\|_0\leq (A\tau, \tau)
$$
for some positive constant $\alpha$. This estimate degenerates
$\alpha\rightarrow 0$ when $\lambda\rightarrow +\infty$.
However the estimate remains true with $\alpha>0$  depending only on $\Omega$
 and $\mu$ if we restrict $\tau$ to functions for which $\d \tau=0$ and
 $\int_{\Omega}\text{tr}(\tau) d\b x=0$, see \cite{Brezzi-Fortin}, also \cite{Arnold-Douglas-Gupta,XieXu2011} for more details.
 \end{remark}

\section{Numerical tests}\lab{s-numerical}

We compute a 2D    pure displacement problem  on the unit square $\Omega=[0,1]^2$
     with a homogeneous boundary condition
         that $u\equiv 0$ on $\partial\Omega$.
In the computation, we let $\mu=1/2$ and $\lambda=1$, and  the  exact solution be
   \e{\lab{e-2}   u= \p{ e^{x-y}x (1-x) y (1-y)\\
    \sin(\pi x)\sin(\pi y)}. }
The true stress function $\sigma$
     and the load function $f$ are defined by the equations in
    \meqref{eqn1},   for the given  solution $u$.

In the computation, the level one grid consists of two right triangles,
   obtained by cutting the unit
   square with a north-east line.
Each grid is refined into a half-sized grid uniformly, to get
   a higher level grid.
In all the computation, the discrete systems of equations are
solved by Matlab backslash solver.

\begin{table}[htb]
  \caption{  The errors, $\epsilon_h =\sigma -\sigma_h$,
     and the order of convergence, by the $P_3$ element, for \meqref{e-2}.}
\lab{b-2}
\begin{center}  \begin{tabular}{c|cc|cc|cc|cc}  
\hline & $ \|u-u_h\|_{0}$ &$h^n$ &
    $ \|\epsilon_h\|_0$ & $h^n$  &
    \multispan{1} $ \|\d\epsilon_h \|_0$  & $h^n$ &
      \multispan{1} $\dim V_h$ &$\dim \Sigma_h$   \\ \hline
 1& 0.118116&0.00& 0.89740816&0.00& 4.917949&0.00&    24&     50\\
 2& 0.024156&2.29& 0.14324287&2.65& 0.981834&2.32&    96&    163\\
 3& 0.002462&3.29& 0.01069158&3.74& 0.132268&2.89&   384&    587\\
 4& 0.000285&3.11& 0.00069804&3.94& 0.016842&2.97&  1536&   2227\\
 5& 0.000035&3.03& 0.00004416&3.98& 0.002115&2.99&  6144&   8675\\
      \hline
\end{tabular}\end{center} \end{table}

First, we use the $P_3$ finite element, $k=3$ in \meqref{Vh} and \meqref{Sh},
   i.e., the $P_3$ stress element and $P_2$ displacement element.
In Table \mref{b-2}, the errors and the convergence order
   in various norms are listed  for the true solution \meqref{e-2}.
An order 3 convergence is observed  for both displacement and stress,
  see Table \mref{b-2}, as shown in the theorem.
For better observing this property,  we plot the
   finite element solution $(\sigma_h)_{11}$ and its error,
    on level 4 grid, in Figure \mref{f-1}.
We also plot   the
   finite element solution $(u_h)_{1}$ and its error,
    on level 4 grid, in Figure \mref{f-2}.
It is apparent that there is at least one superconvergent point on
   each triangle, for the $P_3$ solutions,  but not for the $P_4$ solutions.

  \begin{figure}[htb]\setlength\unitlength{1in}\begin{center}
    \begin{picture}(3.8,4.5)(.1,.1)
  \put(0,2){ \epsfysize=1.52in \epsfbox{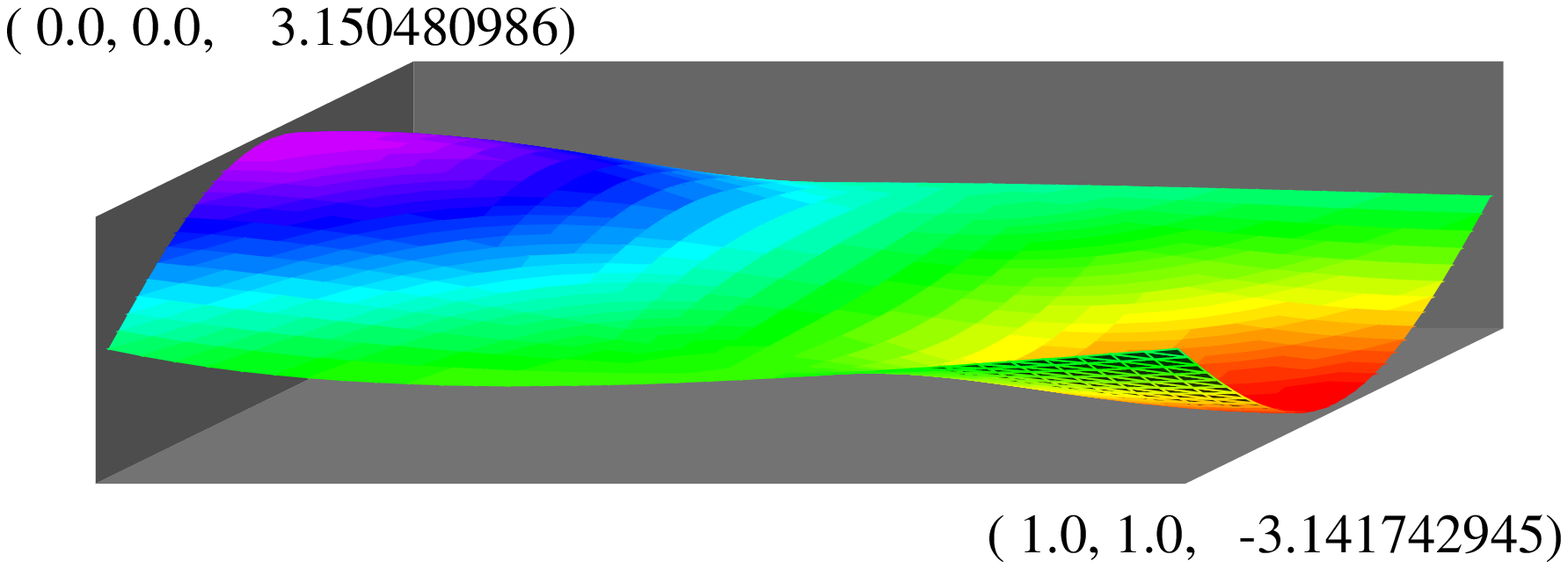}}
  \put(0,0.5){ \epsfysize=1.52in \epsfbox{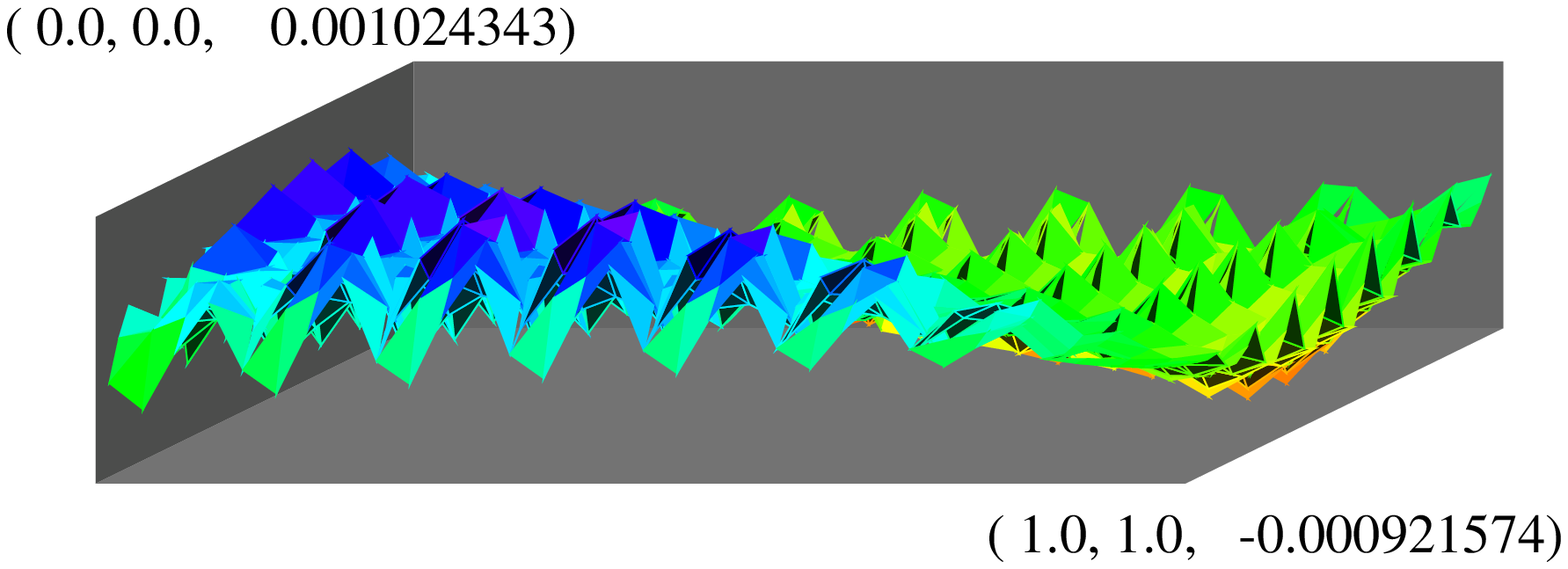}}
  \put(0,-1){ \epsfysize=1.52in \epsfbox{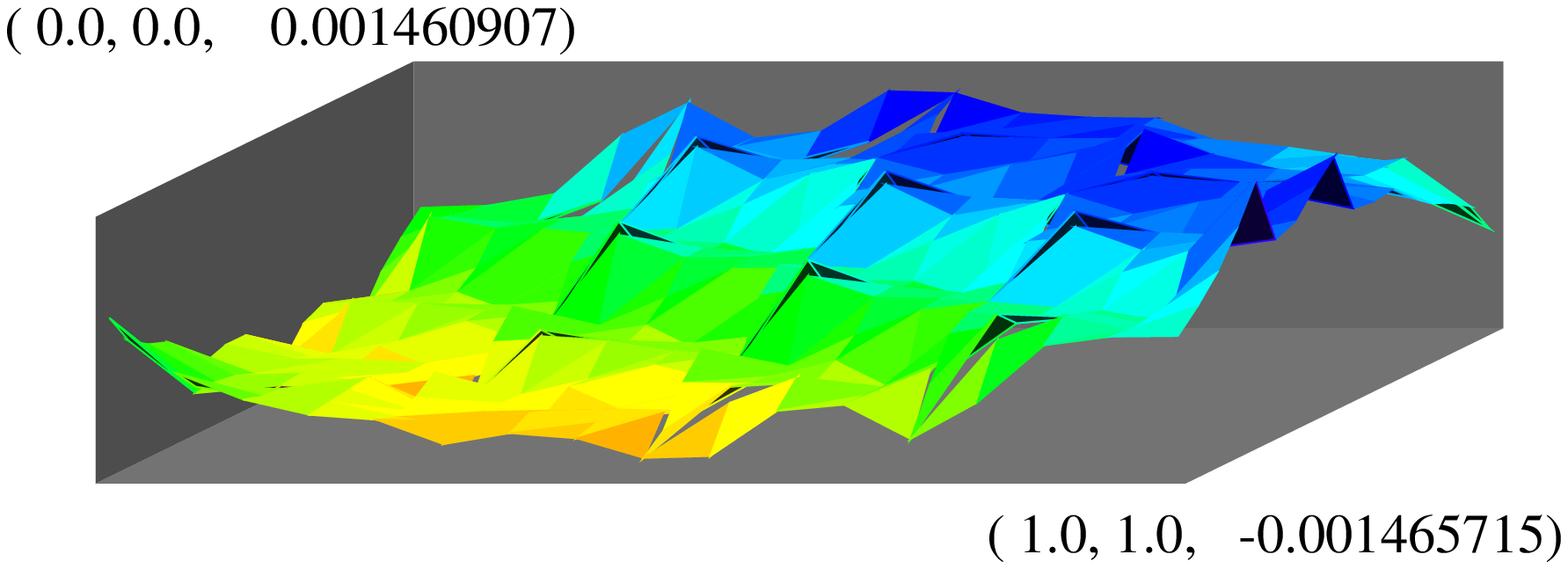}}
    \end{picture}
\caption{The solution of $(\sigma_h)_{11}$ and the error by $P_3$ finite element on
      level 4.  The error (bottom) for $(\sigma_h)_{11}$  by $P_4$ finite element on
      level 3. }
\lab{f-1}\end{center}
\end{figure}

In the second computation, we use the $P_4$ finite element,  i.e.,
$k=4$ in \meqref{Vh} and  \meqref{Sh}.
The data are listed in Table \mref{b-3}.
This time,  the order of convergence is exactly as proved in the theorem,
  order 4 in all norms. It can  be seen from  \mref{f-2}
   that $P_4$ solutions have no zero point (superconvergent point) for $u$ on each
   element.

\begin{table}[htb]
  \caption{   The errors, $\epsilon_h =\sigma -\sigma_h$,
     and the order of convergence, by the $P_4$ element  ($k=4$ in \meqref{Vh} and
    \meqref{Sh}), for \meqref{e-2}.}
\lab{b-3}
\begin{center}  \begin{tabular}{c|cc|cc|cc|cc}  
\hline & $ \|u-u_h\|_{0}$ &$h^n$ &
    $ \|\epsilon_h\|_0$ & $h^n$  &
    \multispan{1} $ \|\d\epsilon_h \|_0$  & $h^n$ &
      \multispan{1} $\dim V_h$ &$\dim \Sigma_h$   \\ \hline
   1& 0.04847978&0.0& 0.162921&0.0& 1.454469&0.0&    40&     78\\
 2& 0.00288821&4.1& 0.005690&4.8& 0.085544&4.1&   160&    267\\
 3& 0.00019094&3.9& 0.000199&4.8& 0.005586&3.9&   640&    987\\
 4& 0.00001211&4.0& 0.000007&4.9& 0.000353&4.0&  2560&   3795\\
      \hline
\end{tabular}\end{center} \end{table}

  \begin{figure}[htb]\setlength\unitlength{1in}\begin{center}
    \begin{picture}(3.8,4.5)(.1,.1)
  \put(0,2){ \epsfysize=1.52in \epsfbox{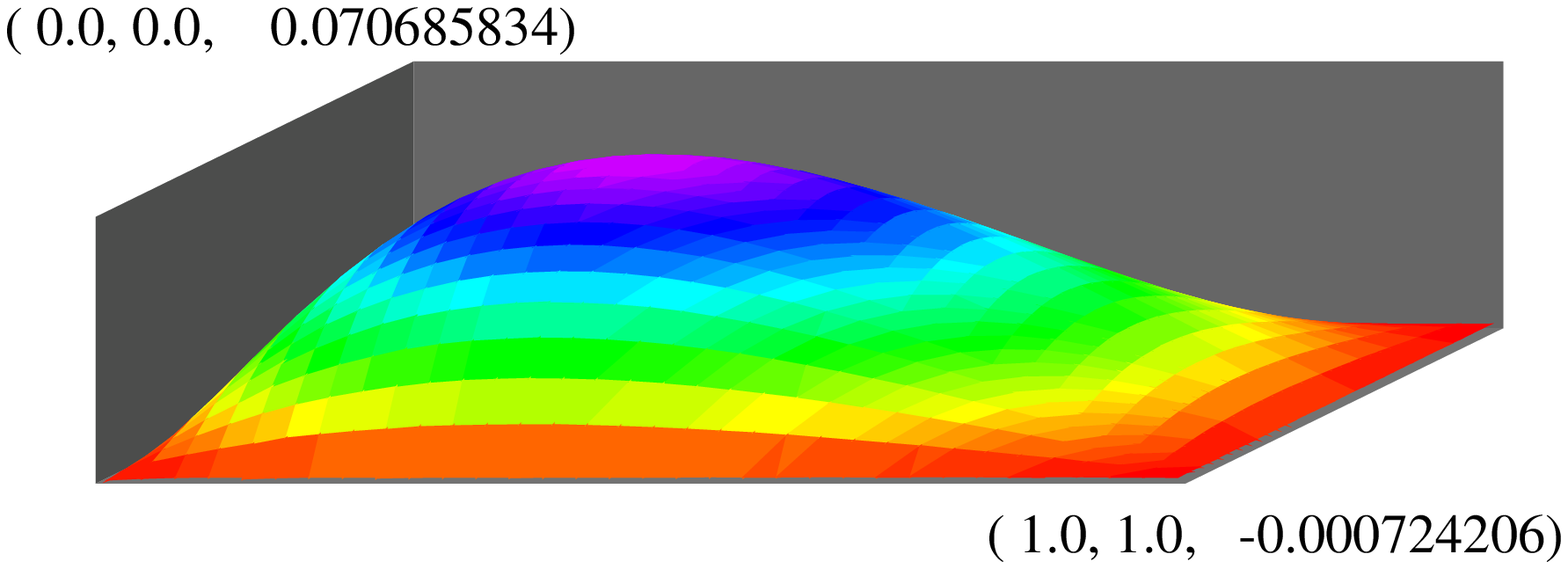}}
  \put(0,0.5){ \epsfysize=1.52in \epsfbox{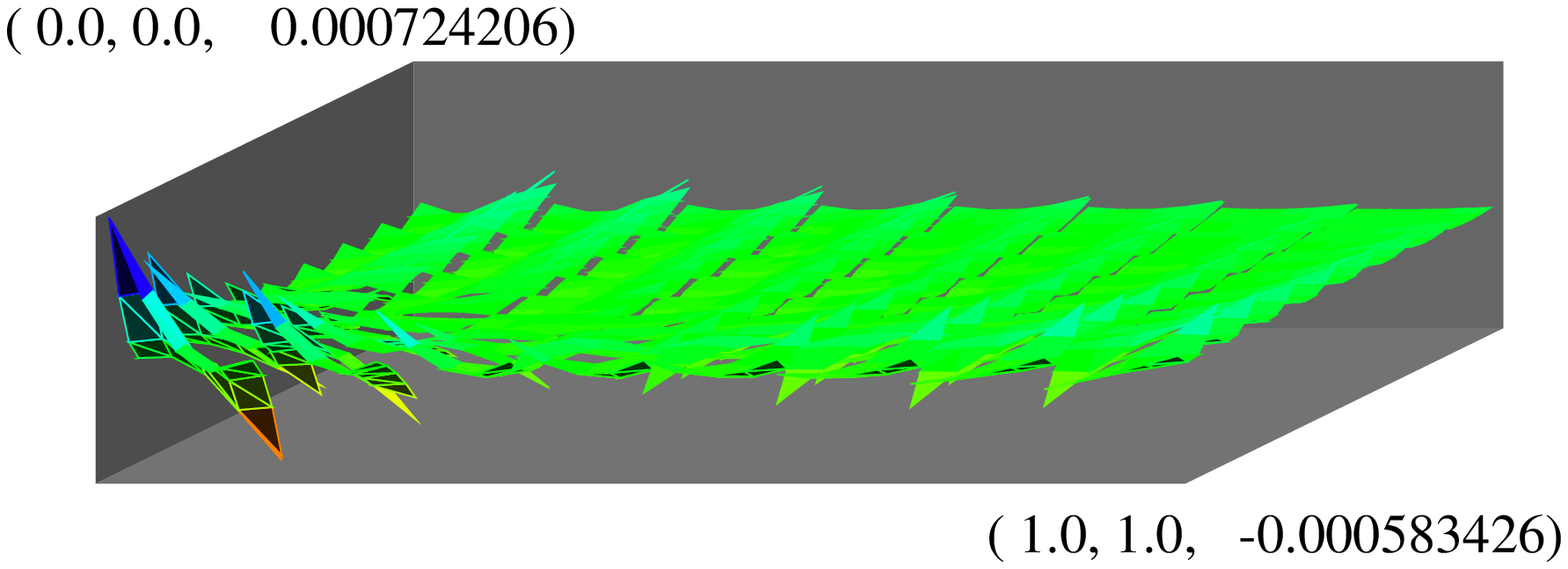}}
  \put(0,-1){ \epsfysize=1.52in \epsfbox{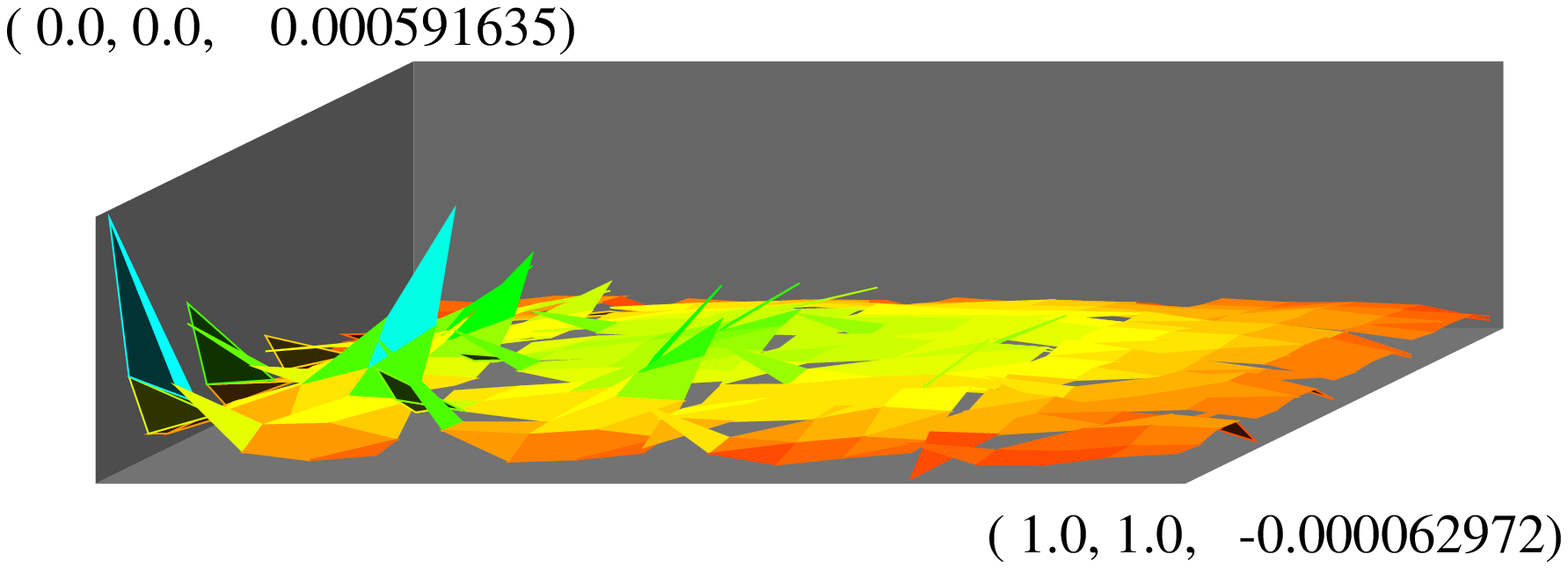}}
    \end{picture}
\caption{The solution of $(u_h)_{1}$ and the error by $P_3$ finite element on
      level 4.  The error (bottom) for $(u_h)_{1}$  by $P_4$ finite element on
      level 3. }
\lab{f-2}\end{center}
\end{figure}

In the third computation, we use the $P_5$ finite element,  i.e.,
$k=5$ in \meqref{Vh} and  \meqref{Sh}.
The data are listed in Table \mref{b-4}.
The order of convergence in $H(\d)$ norm is as proved in the theorem,
  order 5.
Again,  like the $P_3$ and $P_4$ elements,  the $P_5$ element has a  sixth order convergence
   in $L^2$ for the stress.

\begin{table}[htb]
  \caption{ The errors, $\epsilon_h =\sigma -\sigma_h$,
     and the order of convergence, by the $P_5$ element  ($k=5$ in \meqref{Vh} and
    \meqref{Sh}), for \meqref{e-2}.}
\lab{b-4}
\begin{center}  \begin{tabular}{c|rr|rr|rr|rr}  
\hline & $ \|u-u_h\|_{0}$ &$h^n$ &
    $ \|\epsilon_h\|_0$ & $h^n$  &
    \multispan{1} $ \|\d\epsilon_h \|_0$  & $h^n$ &
      \multispan{1} $\dim V_h$ &$\dim \Sigma_h$   \\ \hline
 1& 0.0053888&0.0& 0.022720&0.0& 0.243435&0.0&    60&    112\\
 2& 0.0005013&3.4& 0.002159&3.4& 0.019784&3.6&   240&    395\\
 3& 0.0000145&5.1& 0.000040&5.7& 0.000655&4.9&   960&   1483\\
 4& 0.0000004&5.0& 0.000001&5.9& 0.000021&5.0&  3840&   5747\\      \hline
\end{tabular}\end{center} \end{table}

In the last computation,  we use the $P_3$ Arnold--Winther element
    \cite{Arnold-Winther-conforming} where
   the stress space is the $P_3$ polynomials whose divergence is $P_1$.
      The total degrees of freedom for the stress for the new $P_3$ element are $3|\mathbb{V}|+4|\mathbb{E}|+9|\mathbb{K}|$, where    $|\mathbb{V}|$, $|\mathbb{E}|$, and $|\mathbb{K}|$ are the numbers of vertices, edges and elements of $\mathcal{T}_h$, respectively, while those for the Arnold--Winther element are $3|\mathbb{V}|+4|\mathbb{E}|+3|\mathbb{K}|$.  Since the
     nine bubble functions on each element can be easily condensed,  these two elements almost have  the same complexity for solving. Nevertheless, the new element has  one order higher convergence than the Arnold-Winther element, see  the data in Tables \mref{b-2} and \mref{b-5}.

\begin{table}[htb]
  \caption{ The errors, $\epsilon_h =\sigma -\sigma_h$,
     and the order of convergence, by the $P_3$ Arnold-Winther element
      \cite{Arnold-Winther-conforming} , for \meqref{e-2}.}
\lab{b-5}
\begin{center}  \begin{tabular}{c|rr|rr|rr|rr}  
\hline & $ \|u-u_h\|_{0}$ &$h^n$ &
    $ \|\epsilon_h\|_0$ & $h^n$  &
    \multispan{1} $ \|\d\epsilon_h \|_0$  & $h^n$ &
      \multispan{1} $\dim V_h$ &$\dim \Sigma_h$   \\ \hline
1&  0.27384&0.0&  1.21549&0.0&   6.97007&0.0&   12&    38\\
2&  0.07429&1.9&  0.16642&2.9&   2.13781&1.7&   48&   115\\
3&  0.01959&1.9&  0.02180&2.9&   0.57734&1.9&  192&   395\\
4&  0.00497&2.0&  0.00274&3.0&   0.14709&2.0&  768&  1459\\
5&  0.00125&2.0&  0.00034&3.0&   0.03694&2.0& 3072&  5603\\      \hline
\end{tabular}\end{center} \end{table}

\end{document}